\input epsf
\input amstex
\magnification 1200
\documentstyle{amsppt}
\NoBlackBoxes
\def\1{\bold 1}
\def\a{\alpha}
\def\l{\lambda}
\def\g{\frak g}
\def\h{\frak h}
\def\Z{\Bbb Z}
\def\C{\Bbb C}
\def\Q{\Bbb Q}
\def\R{\Bbb R}

\def\d{\partial}
\def\i{\text{i}}

\def\Id{\operatorname{Id}}

\def\dim{\operatorname{dim }}
\def\Hom{\operatorname{Hom}}

\def\v{^\vee}
\def\<{\langle}
\def\>{\rangle}
\def\[{[\![}
\def\]{]\!]}

\def\sltwo{{\frak s\frak l _2 }}
\def\sltwohat{\widehat{\frak s\frak l _2 }}
\def\sln{\frak{sl}_n}
\def\Ug{U_q\frak g}

\def\o{\otimes}
\def\odot{\operatornamewithlimits{\overset \bold{.} \to \otimes}}

\def\ghat{{\hat{\frak g}}}
\def\gtilde{{\tilde{\frak g}}}
\def\Ok{{\Cal O^{int}_\varkappa}}

\topmatter
\title On inner product in modular tensor categories. II\\
 Inner product on conformal blocks and affine inner product identities.
\endtitle

\author 
Alexander A. Kirillov, Jr.\\
\it Dept. of Mathematics, MIT\\
Cambridge, MA 02139, USA\\
 e-mail: kirillov\@math.mit.edu\\
http://web.mit.edu/kirillov/www/home.html
\endauthor
\leftheadtext{Alexander Kirillov, Jr.}
\rightheadtext{On inner product on conformal blocks}
\date November 5, 1996\enddate

\endtopmatter

\document

\head Introduction. \endhead

This is the second part  of the paper \cite{K}, and we freely use
the notations from there. In \cite{K}, we defined for every modular
tensor category (MTC)  inner products on the spaces of 
morphisms and proved that the inner product on the space $\Hom
(\bigoplus X_i\o X^*_i, U)$  is modular invariant.  Also, we have
shown that in the case of the MTC arising from the representations of
the quantum group $U_q \sln$ at roots of unity and $U$ being a
symmetric power of the fundamental representation, the action of
modular group on the space of intertwiners $\Hom
(\bigoplus X_i\o X^*_i, U)$ can be written explicitly in terms of
Macdonald's polynomials at roots of unity.

In this paper, we apply the same construction to the MTC $\Ok$ coming
from the integrable representations of affine Lie algebras or,
equivalently, from Wess-Zumino-Witten model of conformal field
theory. In Sections~8, 9 we briefly recall the construction of this
category, first suggested by Moore and Seiberg (see \cite{MS1,2}) and
later refined by Kazhdan and Lusztig (\cite{KL1--4}) and Finkelberg
(\cite{F}). In particular, spaces of morphisms in this category are
the spaces of conformal blocks of the WZW model. We prove in Section
10 that $\Ok$ is hermitian, i.e. can be endowed with a suitable
complex conjugation. Thus, the general theory developed in Section~2
of \cite{K} gives us an inner product on the spaces of conformal
blocks, and so defined inner product is modular invariant. We show in
Section 11 that this definition of the inner product is constructive:
it can be rewritten so that it only involves Drinfeld associator, or,
equivalently, asymptotics of solutions of the  Knizhnik-Zamolodchikov
equations.

Since there are integral formulas for the solutions of the KZ equations
(\cite{SV}), this shows that the inner product on the space of
conformal blocks can be written  in terms of asymptotics of
certain integrals. In the case $\g=\sltwo$ these asymptotics can be
calculated (see \cite{V1}), using Selberg integral, and the answer is
given by certain products of $\Gamma$-functions. Thus, in this case we
can write explicit formulas for the inner product on the space of
conformal blocks; we do it in Section~12.  These expressions are
closely related with those suggested by Gawedzki et al. (\cite{G,
FGK}), though their approach is completely different from ours. There
is little doubt that the same correlation holds for arbitrary Lie
algebras. 

The  norms of intertwiners for $\sltwohat$ are given by a formula
which is very much similar to the Macdonald inner product formula,
which gives the norms of intertwining operators for
$U_q\sltwo$ (see \cite{K}). Motivated by this, we formulate a
conjecture about the norms of the intertwining operators (in some
special cases, which are related with  Macdonald's theory) for
$\widehat{\sln}$. 

 It was proved recently in \cite{KL1--4, F} (though this statement was
widely believed long before) that as a braided category $\Ok$ 
is equivalent to the reduced category of representations of $\Ug$ with
$q$ being a root of unity, which was discussed in the first part of
this paper. Since in both categories the action of modular group is
written in terms of braiding, they are also equivalent as modular
tensor categories. In particular, this implies  that if we let
$\g=\sln$ and consider the space of conformal blocks on the
torus with one puncture, to which a symmetric power of fundamental
representation is assigned, then in some basis the action of the
modular group  on this space is 
given by special values of Macdonald's polynomials. Thus, for
$\g=\sltwo$ the action of modular group on any space of conformal
blocks on a torus  can be written in some basis in terms of
$q$-ultraspherical polynomials, and this basis is related to the
standard one by a matrix of gamma-functions.

\head 8. Integrable representations of affine Lie algebras and conformal
blocks \endhead

Here we briefly recall the main definitions, referring the reader to
\cite{MS1, TUY, KL} for details. This section is completely expository.

\subhead Integrable modules \endsubhead

As before, let $\g$ be a simple Lie algebra over $\C$ with the
invariant bilinear form $(\, , \,)$ normalized so that $(\theta,
\theta)=2$ where $\theta$ is the highest root. Let $\gtilde=\g\o \C[t,
t^{-1}]\oplus\C c$ be the corresponding affine Lie algebra; as usual,
we denote $x[n]=x\o t^n$. We also denote by $\ghat$ its completion:
$\ghat=\g\o\C((t))\oplus\C c$.  In this whole paper, we fix a
non-negative integer $k$ (level). Let $\varkappa=k+h\v$, and let $\Ok$
be the category of finite-length integrable highest-weight
$\ghat$-modules of level $k$. This category is semisimple and the
simple objects in this category are the irreducible modules $L_{\l,
k}, \l\in C$, where

$$C=\{\l\in P^+|\<\l, \theta\v\>\le k\}=\{\l\in P^+|\<\l+\rho,
\theta\v\><\varkappa\},\tag 8.1$$
(see \cite{Kac}). 

There exists a natural notion of duality in this
category: for $V\in \Ok$, let $DV$ be the restricted  dual space to
$V$ (that is, $DV$ is direct sum of spaces dual to homogeneous
components of $V$) and the action of $\gtilde$ is defined as the usual
action in the dual space twisted by the automorphism
$\sharp$ defined by  
$$x[n]^\sharp=(-1)^n x[-n], \quad c^\sharp=-c.$$ 

 It is easy to check that $DV\in \Ok$, and $D(DV)$ is canonically
isomorphic to $V$: the usual isomorphism of vector spaces
$V^{**}\simeq V$ is $\ghat$-isomorphism.  Also, $DL_{\l, k}\simeq
L_{\l^*, k}$, though the isomorphism is not canonical.

\subhead  Conformal blocks\endsubhead

Let $\Cal X$ denote the following collection of data:
\roster
\item $X$ -- a non-singular compact complex curve, 

\item $z_1, \dots, z_n$ -- distinct point on $X$ divided into two sets
$In$ and $Out$,

\item  $w_i$ -- local parameter near the point $z_i$. i.e. a
holomorphic function in a neighborhood of $z_i$ such that $w_i(z_i)=0,
w'_i(z_i)\ne 0$. 
\endroster

 We assume that on each connected component of $X$ there is at least
one of the points $z_i$.

With each point $z_i$ we associate a Lie algebra $\g\o\C((w_i))$. Let

$$\gathered 
\ghat_{in}=\left(\bigoplus_{i\in In}\g\o\C((w_i))\right)\oplus \C c\\
\ghat_{out}=\left(\bigoplus_{i\in Out}\g\o\C((w_i))\right)\oplus \C c.\\
\endgathered
\tag 8.2$$
and the cocycle defining the central extension in $\ghat_{in}$ 
(respectively, $\ghat_{out}$)
is the sum of standard cocycles on each of $\g\o \C((w_i)), i \in  In$
(respectively, $i\in Out$).

Also, let us consider the Lie algebra $\Gamma$:
$$\Gamma(\Cal X)=\{\g-\text{valued meromorphic functions on $X$ regular
outside of } z_1,\dots, z_n\}\tag 8.3$$
and its central extension $\widehat\Gamma=\Gamma\oplus \C c$ with the
defining cocycle given by 
$$c(f,g)=\sum_{i\in In}\text{Res}_{z_i} g\,df.$$

Expanding a function $f\in \Gamma$ near a point $z_i$ in a
Laurent series in $w_i$ we get a Lie algebra homomorphism 
$\pi_i:\Gamma(\Cal X)\to \g\o \C((w_i))$ (if $\Cal X$ is connected,
this is an embedding). 

Taking direct sum over all $i\in In$ (respectively, $i\in Out$), we
get embeddings 

$$\gathered
\pi_{in}: \widehat \Gamma(\Cal X)\subset \ghat_{in}:\qquad
f\mapsto \oplus_{i\in In} \pi_i(f),\qquad c\mapsto c,\\
\pi_{out}: \widehat \Gamma(\Cal X)\subset \ghat_{out}:\qquad
f\mapsto \oplus_{i\in Out} \pi_i(f),\qquad c\mapsto -c.
\endgathered\tag 8.4$$

One can easily check that these embeddings are Lie algebra
homomorphisms.

\definition{Definition 8.1} Let $\Cal X$ be as above, and assume that
we are given integrable modules $V_1,\dots, V_n \in \Ok$ assigned to
the points $z_1, \dots, z_n$ respectively. Let us consider $V_i$ as a
module over $\g\o \C((w_i))\oplus \C c$. Then the corresponding space
of conformal blocks is defined by

$$\aligned W(\Cal X; V_1, \dots, V_n)=
\{\Phi: 
	\o_{i\in In}V_i\to \widehat\o_{i\in Out} V_i
	\bigl|&(\pi_{out}(f))^\sharp\Phi=\Phi\pi_{in}(f)\\
	&\text{ for all }f\in \Gamma(\Cal X)
\}\endaligned
\tag 8.5$$
where $\widehat \o$ is the completion of the tensor product with
respect to the homogeneous grading.
\enddefinition

\remark{Remark 8.2} If $i\in Out$ then denote by $\Cal X'$ the same data
as $\Cal X$ except that now we consider $i$ as an element of $In$:
$In'=In\cup \{i\}, Out'=Out\setminus\{i\}$. Then we have a canonical
isomorphism: 

$$W(\Cal X; V_1, \dots, V_n)=W(\Cal X'; V_1, \dots,DV_i,\dots V_n).$$

Thus, it suffices to consider conformal blocks when all points are
incoming (or all are outgoing), as is done in \cite{TUY}; however, it
is more convenient to consider the general situation. 

In a similar way, it can be shown that if $\Cal X'$ is obtained from
$\Cal X$ by marking one more incoming point $z_0$ and assigning to it the
representation  $V_0=L_{0,k}$ then the map 

$$\aligned 
W(\Cal X'; V_0, V_1, \dots, V_n)&\to W(\Cal X; V_1, \dots, V_n)\\
\Phi&\mapsto \Phi(1,\dots),
\endaligned$$
where $1$ is the highest weight vector in $L_{0,k}$, is an
isomorphism. Thus, the condition that there is at least one point on
each connected component of $X$ is inessential.
\endremark

\remark{Remark 8.3} In fact, the definition of conformal blocks works
in more general situation: we can allow $X$ to be a semi-stable
singular curve. This plays a crucial role in proving the gluing axiom
(\cite{TUY}). However, all of the results we need in this paper can be
formulated without reference to singular curves.
\endremark

\example{Example 8.4} Let $n=3, X=\C P^1$ with global coordinate $w$, 
$In=\{\infty\}, Out=\{0,z\}$ with local parameters $1/w, -w, z-w$
respectively. Let $\Phi: V_\infty \to V_0\hat \o V_z$ be an element of
the space of conformal blocks $ W(\Cal X; V_0, V_z, V_\infty)$.  Then 

$$\Phi x[n] =\biggl(  x[n]\o 1 + 
				\sum_{i\ge 0} { -n \choose i} z^{-n-i} 
				 1\o  x[-i]
			\biggr)\Phi.
$$  

This slightly differs from the usual formulas in the  physical literature
where usually $\infty$ is considered as an outgoing  and $0, z$ as
incoming points.
\endexample 
 
The following result is well known; we refer the reader to \cite{TUY}
for the proof.

\proclaim{Proposition 8.5} The spaces of conformal blocks are always
finite-dimensional. \endproclaim

\subhead Correlation functions \endsubhead

For a simple module $L_{\l,k}$ let $L_{\l,k}[0]$ be the $\g$-module
generated by the highest weight vector; clearly, $L_{\l,k}[0]\simeq
L_\l$. It is easy to show that this operation can be defined in
invariant terms; thus, we have a faithful functor $V\mapsto V[0]$
from  $\Ok$ to the category $Rep\ \g$ of finite-dimensional
representations of $\g$. Note that we have canonical identifications 
 $(DV)[0]\simeq (V[0])^*, L_{0, k}[0]\simeq \C$. Denote by $i:V[0]\to
V$ and $p:V\to V[0]$ the canonical embedding and projection
respectively.

For every $\Phi\in W(\Cal X;V_1,\dots, V_n)$ define its correlation function
(or, which is the same, its highest term) $\<\Phi\>$ by

$$\<\Phi\>=p\o \dots\o p\circ \Phi\circ (i\o \dots\o i)
	\in \Hom_\g(\o_{In} V_i[0],\o_{Out} V_i[0]).
\tag 8.6$$

\proclaim{Proposition 8.6} On a sphere, every $\Phi\in  W(\Cal
X,V_1,\dots, V_n)$ is uniquely defined by its
correlation function $\<\Phi\>$. In other words, we have an embedding:

$$  \gathered
 W(\Cal X;V_1,\dots, V_n) \subset
	\Hom_\g(\bigotimes_{In} V_i[0],\bigotimes_{Out} V_i[0]) \\
\Phi\mapsto \<\Phi\>.\endgathered\tag 8.7$$

\endproclaim

\remark{Remark} This statement is not true on higher genus surfaces.
\endremark 

Note that this embedding depends on the choice of points $z_i$ and
local parameters at these points (in fact, it depends only on 1-jet of
the local parameter). For irreducible  $V_i$ the  image of this
embedding  can be explicitly described; we refer the
reader to \cite{FSV} (for the case of three-punctured sphere, this
description also appears in \cite{TUY}).

\subhead Flat connection \endsubhead

The definition of conformal blocks which we gave 
depends on $\Cal X$; thus, we can consider the space of
conformal blocks as  a finite-dimensional vector bundle over the
corresponding moduli space. This vector bundle is called the bundle of
conformal blocks.  
 
It turns out that this vector bundle  has a natural flat connection,
which was first calculated by Knizhnik and Zamolodchikov on the sphere
and by Bernard on a torus. This
connection can be naturally defined using the Sugawara construction;
we refer the reader to  \cite{TUY, KL2} for details, giving here only the
answer.

From now on we assume the following 

$$
\Cal X=\left(
\gathered
X=\C P^1, w - \text{global coordinate on } X\\
In=\{\infty\}, \quad w_\infty=1/w,\\
Out=\{z_1,\dots, z_n\}, \quad w_i=z_i-w,
\endgathered\right) 
\tag  8.8$$

By Proposition~8.6, for each $z_1, \dots, z_n$ the space $W(\Cal
X,V_1, \dots,V_n, V_\infty) $ can be identified with a certain
subspace in $\Hom_\g(V_\infty[0], V_1[0]\o \dots\o V_n[0])$. Thus, to
define a flat connection on the space of conformal blocks it suffices
to define a flat connection in the trivial vector bundle with the fiber
$\Hom_\g(V_\infty[0], V_1[0]\o \dots\o V_n[0])\simeq (V_1[0]\o \dots\o
V_n[0]\o V^*_\infty[0])^\g$ which would preserve the subbundle of conformal
blocks. Such a connection is obtained from the Knizhnik-Zamolodchikov
connection on $V=V_1[0]\o \dots\o V_n[0]$

$$(k+h\v)\frac{\d}{\d z_i} \psi =\left(\sum \Sb j=1\dots n\\j\ne
i\endSb\frac{\Omega_{ij}}{z_i-z_j} \right)\psi,
\tag 8.9$$ 
where $\Omega$ is the standard $\g$-invariant element in $\g\o
\g$: if $x_i, x^i$ are dual bases in $\g$ with respect to the inner
product $(\, , \, )$ then $\Omega=\sum x_i\o x^i$. As usual, we use
the notation $\Omega_{ij}=\pi_i\o \pi_j(\Omega), i,j=1,\dots, n$.

This connection can be extended to a connection with values in $V\o
V^*_\infty[0]$ with trivial action on the last factor. Furthermore,
since this connection commutes with the action of $\g$, it also
defines a connection on the trivial vector bundle with the fiber
$(V_1[0]\o \dots\o V_n[0]\o V^*_\infty[0])^\g$, and it can be checked
that it preserves the subbundle of conformal blocks (see
\cite{FSV1,2}).

\example{Example 8.7} Let us consider the conformal blocks on a
3-punctured sphere: $In=\{\infty\},Out=\{z_1, z_2\}$ and 
 assume that the modules $V_i$
are irreducible: $V_\infty=L_{\l,k}, V_1=L_{\mu,k},V_2=L_{\nu,k}$. In
this case it is easy to check that a section $\Phi(z_1, z_2)$ of the
bundle of conformal blocks  is flat iff
$\<\Phi(z_1, z_2)\>=(z_1-z_2)^{\Delta_\l-\Delta_\mu-\Delta_\nu} g$, where
$\Delta_\l=\frac{(\l, \l+2\rho)}{2(k+h\v)}$, and $g\in
\Hom_\g(L_\l,L_\nu \o L_\mu)$ does not depend on $z_i$. 

 In this case we will use the notation $\Phi^g(z), z=z_1-z_2$ for such
a flat section; the operators $\Phi(z)$ are called chiral vertex
operators.
 
\endexample

\head 9. Category $\Ok$ as modular category \endhead

In this section we recall the construction of a tensor (and in fact,
modular tensor) structure on the category $\Ok$. This section is again
expository; we refer the reader to \cite{MS, KL} for details and
proofs.

 As before, let us assume that $\Cal X$ is $n+1$-punctured sphere
(8.8) and $n\ge 1$. The moduli space of all such punctured spheres is
the configuration space 

$$X_{n}=\{(z_1,\dots, z_n)\in \C^n|z_i\ne z_j\}.$$

As was discussed above, the space of conformal blocks is a local
system over this space with the connection given by the KZ equations
(8.9). Note also that these equations imply that every flat section of
this local system is invariant under translations $(z_1, \dots,
z_n)\mapsto (z_1+c, \dots, z_n+c)$ and thus consideration of this
local system is can be reduced to consideration of a local system on
$X^0_n=\{\bold z\in X_n|z_1=0\}$. Define

$$\Cal D_n=\{(z_2,\dots, z_n)\in \C^{n-1}|0<|z_2|<|z_3|<\dots
<|z_n|,z_i\notin \R_{<0}\}\subset X_n,
\tag  9.1$$

Note that $\Cal D$ is contractible.

Now we can formulate the main theorem of this section, which is
essentially due to Moore and Seiberg (see also \cite{KL}).

\proclaim{Theorem 9.1} The category $\Ok$ can be endowed with the 
structure of a ribbon tensor category such that:

\roster
\item If we denote the tensor product in this category by $\odot$
\rom{(}to avoid confusion with the usual tensor product of vector
spaces\rom{)} then for any $n\ge 1$ and for any choice of
representations $V_1, \dots, V_n, V_\infty\in \Ok$

$$\Hom_\Ok(V_\infty, V_1\odot \dots\odot V_n)=\Gamma(\Cal D_n, W),
\tag 9.2
$$ 
where $W=W(\Cal X; V_1, \dots, V_\infty)$ is the local system of
conformal blocks on $n+1$-punctured sphere \rom{(8.8)}, and $\Gamma$
stands for the space of global sections of this local system over
$\Cal D_n$ defined by \rom{(9.1)}. 

\item The unit object is $\1=L_{0,k}$ and the maps $V\simeq \1\odot V$
are constructed as in Remark~\rom{8.2}.

\item The dual object  is given by
$V^*=DV$, and  the maps $i_V:\1\to V\odot V^*$ are defined so that 
$i_{L_{\l,k}}= \Phi^{i_\l}$, where $\Phi^g$ is defined in
Example~\rom{8.7}, $i_\l$ is the canonical map of $\g$-modules
$\C\to L_\l\o L_\l^*$.

\item The  isomorphism $\delta_V: V\to V^{**}=D(DV)$ coincides with
canonical identification of vector spaces $V\simeq V^{**}$
\rom{(}recall that as a vector space $DV=V^*$ is the restricted dual
to $V$\rom{)}, and the twist $\theta$ is given by $\theta=e^{2\pi\i
L_0}$, where $L_0$ is the Sugawara element in the completion of
$U\ghat$; thus,

$$\theta|_{L_{\l,k}}=e^{2\pi\i \Delta_\l},\qquad 
	\Delta_\l=\frac{(\l, \l+2\rho)}{2(k+h\v)}.
\tag 9.3$$

\endroster
\endproclaim

\remark{Remark} 
 Since $\Cal D_n$ is simply connected, we could just say that
we fix some  particular choice of points $z_1, \dots, z_n$, say,
$\bold z=(0,1,\dots, n-1)$ and let 
$\Hom_\Ok(V_\infty, V_1\o \dots\o V_n)=W(\Cal X; V_1, \dots,
V_\infty)$. 
\endremark

We do not define here the associativity and commutativity
isomorphisms, referring the reader to the original papers. However,
it is necessary to mention that the construction of the associativity
isomorphism

$$\Hom_\Ok(L_{\l,k}, (V_1\odot V_2)\odot V_3) \simeq 
\Hom_\Ok(L_{\l,k}, V_1\odot (V_2\odot V_3))\tag 9.4
$$
is based on the fact that we can identify each of these spaces with
the space of conformal blocks on a 4-punctured sphere (see (9.2)). The
identifications are obtained by considering the asymptotics of the
flat sections in different asymptotic zones. Therefore, the
associtivity morphism is written in terms of asymptotics of solutions
of the  KZ equations in 3 variables. This implies the following lemma.

\proclaim{Lemma 9.2} If $V_1, V_2, V_3\in \Ok, \l\in P_+$ are such that on the
 space of singular  vectors of weight $\l$ in $V_1[0]\o V_2[0]\o
 V_3[0]$ the operators $\Omega_{12}, \Omega_{23}$ commute then the
 associativity isomorphism \rom{(9.4)} is trivial, i.e. coincides with
 the restriction of the associativity isomorphism for vector spaces $V_1[0],
 V_2[0], V_3[0]$.
\endproclaim

As for any MTC,  we can use the language of ribbon graphs  for representing
morphisms in the category $\Ok$; unfortunately, the associativity morphism
(which is highly non-trivial) ``does not show'' in the pictures, so
one must be careful when performing calculations (see \cite{BN} for an
approach allowing to avoid this difficulty). In particular, we
can define the ``quantum dimension'' of a module $V$ in the same way
as we did for an arbitrary MTC in Section~1; we will denote it by
$\dim_\varkappa V$. Note that it has nothing to do with the usual
dimension of $V$, which is infinite.

\example{Example 9.3}The space of morphisms
$\Hom_\Ok(V_\infty,V_1\odot V_2)$ coincides with the space of chiral
vertex operators (see Examples~8.4, 8.7) and can be identified with a
subspace in the space of $\g$-homomorphisms: it follows from
Proposition~8.6 that we have an embedding

$$\gathered
\Hom_\Ok(V_\infty,V_1\odot V_2)\subset \Hom_\g(V_\infty[0],
V_1[0]\o V_2[0])\\
\Phi^g(z)\mapsto \<\Phi^g(1)\>=g. 
\endgathered
\tag 9.5$$
\endexample

Now we can describe the action of the modular group. Recall the object $H$
defined in Section~1 for any MTC; in our case, it is given by 
$H=\bigoplus_{\l}DL_{\l,k}\odot L_{\l,k}$. The following result
immediately follows from the gluing axiom for conformal blocks.

\proclaim{Lemma 9.4} The space 
$$\Hom_\Ok(H, U)=\bigoplus_{\l\in C} \Hom_\Ok(L_{\l,k}, L_{\l,k}\odot U)
\tag 9.6$$
is isomorphic to the space of conformal blocks on a torus with one
puncture to which the representation $U$ is assigned.
\endproclaim

The following result clarifies the meaning of the  action of $SL_2(\Z)$
introduced in Section~1.

\proclaim{Theorem 9.5}\rom{\cite{MS2}} 
The category $\Ok$ is a
modular category in the sense of Definition~\rom{1.3}, and the action
of $SL_2(\Z)$ on $\Hom_\Ok(H, U)$ defined in Theorem~\rom{1.10}
coincides with the natural geometric action of $SL_2(\Z)$ on the space
of conformal blocks on a torus. 

\endproclaim

This theorem can be proved in the general context of 2-dimensional modular
functor, using Kirby calculus; see \cite{Tu} for this approach. 
In fact, it is known that in any modular tensor category we can define
an action  of the mapping class group of any punctured Riemann
surface on the appropriate space of conformal blocks: this automatically
follows from the possibility to define the action of $SL_2(\Z)$.

\subhead Equivalence of categories $\Ok$ and $\Cal C(\g, \varkappa)$
\endsubhead

The following important  theorem, which was widely believed for several
years, was proved   by M.~Finkelberg in his thesis \cite{F}.

\proclaim{Theorem 9.6} There exist numbers $n(\g)$ such that for $k\ge
n(\g)$ the functor $V\mapsto V[0]$ is an equivalence of modular tensor
categories  $\Ok$ and  the reduced category $\Cal
C(\g,\varkappa)$ of representations of the quantum group $\Ug$ at root of
unity, described in Section~\rom{3}, with $\varkappa=k+h\v$. \endproclaim

Finkelberg's proof is based on the series of papers by Kazhdan and
Lusztig (\cite{KL1--4}).  Since in these papers only the case of
a simply-laced $\g$ is considered, his proof formally works only in
these cases.  However, it turns out that the same arguments work in
non simply-laced case as well if we let $q=e^{\pi\i/m\varkappa}$ (note
$m$ in the denominator!), as we did in Section~3; see \cite{L,
Section~8  and Erratum}.

The restriction $k\ge n(\g)$ mentioned in the theorem appears only
for exceptional root systems (see table in \cite{F}), and for all
$\g$, $n(\g)\le 6$; thus, this theorem is automatically satisfied if
$k\ge 6$.  From now on, we assume that $k$ is chosen to satisfy these
conditions.

In particular, this theorem implies that the quantum dimensions of
objects in both categories coincide: 
$$\dim_\varkappa V=\dim_q V[0],\quad q=e^{\pi\i/m\varkappa}.$$

Another approach to the construction of an equivalence of these
categories was initiated by Schechtman and Varchenko, who showed that
one can identify Verma modules over the quantum group with the
homologies of certain local systems on the configuration space (see
\cite{V2}). However, for rational values of $\varkappa$ this approach
requires use of intersection homology (or, equivalently, perverse
sheaves), which is done in a recent series of papers by Schechtman and
Finkelberg \cite{FS}.

\head 10. Hermitian structure on $\Ok$ and
 inner product on the space of conformal blocks. \endhead
\rightheadtext{Hermitian structure on $\Ok$}

In this section we define a hermitian structure (which is a certain
analogue of the complex conjugation) on the category $\Ok$ and use
it to define an inner product on the spaces of morphisms in this
category.

\subhead Hermitian structure \endsubhead

Recall (see Section~1) that a hermitian structure on a tensor category
is a system of maps $\overline{\phantom{T}}:\Hom(V, W)\to \Hom(V^*,
W^*)$ which satisfies certain compatibility conditions (1.22). When the
category is defined over $\C$ we  assume that  these maps
are $\C$-antilinear. Usually, to define such a structure we first define
a functor $\overline{\phantom{T}}$ which is antiequivalence of
categories (i.e., we have canonical isomorphisms $\overline{V\o
W}\simeq \overline{W}\o \overline{V}$) and then show that we have
isomorphisms $\overline{V}\simeq V^*$. We can formalize this setup in
the following lemma, proof of which is trivial.

\proclaim{Lemma 10.1} Let $\Cal C$ be a semisimple ribbon category
over $\C$.  Let us assume that we have the following data:

\roster\item A  functor   $\omega:\Cal C\to \Cal C$  satisfying
the following conditions: 
\itemitem{$\bullet$} It is antilinear: for every $\Phi\in \Hom_{\Cal
C},\a\in \C$,  we have $(\a\Phi)^\omega=\bar\a \Phi^\omega$;

\itemitem{$\bullet$} We have functorial isomorphisms
$V^{\omega\omega}\simeq V^{**}, 
(V\o W)^\omega\simeq W^\omega\o V^\omega, (V^\omega)^*\simeq
(V^*)^\omega$ and $\1^\omega\simeq \1$ compatible with each
other in the natural way.

\itemitem{$\bullet$} If $\a, R$ are the associativity and
commutativity isomorphisms 
in $\Cal C$ then $\a^\omega=\a^{-1},R^\omega=R^{-1}$, i.e. we have 
the following commutative diagram

$$ \CD
 (V\o W)^\omega @>(R_{V,W})^\omega>>  (W\o V)^\omega\\
      @|                                     @|   \\
  W^\omega\o V^\omega @>R^{-1}_{ V^\omega,W^\omega}>>
					V^\omega\o W^\omega,
\endCD
$$
and similarly for $\a$. 

\itemitem{$\bullet$} If $\theta:V\to V$ is the universal twist in
$\Cal C$ then $\theta^\omega=\theta^{-1}$; similarly, if $i_V:\1\to
V\o V^*$ and $e_V: V^*\o V\to \1$ are the duality maps then
$i_V^\omega= (1\o \delta^{-1})i_{V^{\omega *}}, 
e_V^\omega=e_{V^{\omega *}}(\delta^{-1}\o 1)$.

\item 
Isomorphisms $X_i^\omega\simeq X_i^*$
for all simple objects $X_i$ \rom{(}which, again, are compatible with
$V^{\omega\omega}\simeq V^{**}$\rom{)}. 
\endroster

Then the category $\Cal C$ can be
uniquely endowed with a structure of hermitian ribbon category by
identifying $V^\omega\simeq V^*$ so that it is compatible with all the
structures of a ribbon category.

\endproclaim

Note that these conditions imply $\dim V\in \R$ for every object $V$;
vice versa, if we know that $\dim V\in \R$ then we can omit the condition
$e^\omega=e$ and replace part (2)  by the following condition:

\roster\item"($2'$)" 
For every $i$, we are given non-zero homomorphisms 

$\phi_i:\1 \to  X_i \o X_i^\omega$

such that $\phi^\omega=\phi$ (up to identification
$V^{\omega\omega}\simeq V^{**}\simeq V$). 

\endroster

We have implicitly used this construction when defining the hermitian
structure on the category of representations of a quantum group in
Section~4.  Now we do the same for the category $\Ok$. Let
$\omega:\g\to \g$ be the $q=1$ specialization of the antilinear
involution defined in Section~4; on the generators it is given by

$$\aligned 
\omega: 
	& e_i\mapsto  e_{i\v},\quad i=1,\dots, r\\
	&f_i \mapsto f_{i\v},\quad i=1,\dots, r\\
	&h\mapsto -w_0(h),\quad h\in \h.
\endaligned 
\tag 10.1
$$

We extend it to an involution on $\ghat$ by

$$
\omega(x[n])= (-1)^n (\omega(x))[n],\quad \omega(c)=c.
\tag 10.2
$$

Obviously, $\omega$ is an automorphism of Lie algebras. 
Now, for every $V\in \Ok$ define $V^\omega$ to be the same set  as
$V$ but with the action of $\ghat$ twisted by $\omega$. Similarly, for
every $\ghat$-homomorphism $\Phi:V\to W$ between modules in $\Ok$ let
$\Phi^\omega$ be the same map  but considered as a homomorphism 
$V^\omega\to W^\omega$. Obviously, this operation
preserves the composition of morphisms. One easily checks that 
 we have canonical identifications $L_{0,k}^\omega\simeq L_{0,k}, 
(DV)^\omega\simeq D(V^\omega)$ (since $\omega$ commutes with
the automorphism $\sharp$ used in the definition of the dual -- see
beginning of Section~8). Also,  define the isomorphism
$V^{\omega\omega}\simeq V$  by $v\mapsto Zv$, where $Z$ is
the central element in a completion of $U\g$ which was constructed in
Theorem~7.2 of \cite{K}.

The most difficult part is to prove that
$\omega$ is a tensor functor, i.e. to construct
isomorphisms $(V\odot W)^\omega= W^\omega\odot V^\omega$. To
do it, let us return to definition of conformal blocks on Riemann
surfaces. 

Let $\Cal X$ be a Riemann sphere with $n$ marked points $z_1,\dots,
z_n$ (see (8.8)) and representations $V_1,\dots, V_n, V_\infty\in \Ok$
assigned to these points, and let $\Phi: V_\infty\to
V_1\widehat\o\dots \widehat\o V_n$
be a conformal block, i.e. a linear map satisfying commutation
relations (8.5). 

\proclaim{Theorem 10.2} Let $\Cal X, \Phi$ be as above and denote by
$\Phi^\omega$ the same $\Phi$ considered as a map $V_\infty^\omega\to
 V_1^\omega\widehat\o\dots\widehat\o V_n^\omega$. Then $\Phi^\omega$
 is a conformal block on the Riemann surface $\Cal X^\omega=\C P^1$
 with marked points $z'_1=-\bar z_1,\dots, z'_n=-\bar z_n, \infty$,
 local parameters $w_i=z'_i-w, w_\infty=1/w$ at $z'_1,\dots, z'_n,
 \infty$ and representations $V_1^\omega, \dots,
 V_n^\omega,V_\infty^\omega $ assigned to these points.
\endproclaim

\demo{Proof}
Follows from the fact that the map $f(w)\mapsto
\omega(f(-\overline{w}))$ is an isomorphism of the algebras of	
rational functions $\Gamma(\Cal X)\simeq \Gamma(\Cal X^\omega)$ which
were defined by (8.3). 
\qed\enddemo

Recalling that we have identified the space of morphisms $V_\infty\to
V_1\odot V_2$ in $\Ok$ with the space of 
conformal blocks on the sphere with marked points $0, 1,\infty$ (see
Theorem~9.1 and remark after it), we immediately deduce from
Theorem~10.2 the following result:

\proclaim{Theorem~10.3} The map $\Phi\mapsto \Phi^\omega$ defined in
Theorem~\rom{10.2} gives rise to a functorial antilinear isomorphism

$$\Hom_\Ok(V_\infty, V_1\odot V_2)\simeq \Hom_\Ok
(V^\omega_\infty,V_2^\omega\odot V_1^\omega) $$
and thus gives rise to an isomorphism $(V_1\odot V_2)^\omega\simeq 
V_2^\omega\odot V_1^\omega$.

\endproclaim

\example{Example 10.4} Let $\Phi=\Phi^g\in \Hom_\Ok(V_\infty, V_1\odot
V_2)$ (see Example~8.7). Then $(\Phi^g)^\omega=\Phi^{(g^\omega)}$,
where 

$$g^\omega:(V_\infty[0])^\omega\to
	V_2[0]^\omega\o V_1[0]^\omega$$
is defined using the involution $\omega$ on $\g$. 
 \endexample

Finally, let us fix for every $\l$ a homomorphism of $\g$-modules
$\varphi_\l: \C\to L_\l\o L_\l^\omega$ such that it is symmetric:
$\varphi_\l(1)^\omega =\varphi_\l (1)$. This defines $\varphi_\l$
uniquely up to a real constant.  Define $\Ok$-morphisms $\phi_\l:
\1\to L_{\l,k}\odot L_{\l, k}^\omega$ by $\phi_\l=\Phi^{\varphi_\l}$.

\proclaim{Proposition 10.5} The functor $\omega$ and the system of
maps $\phi_\l: \1\to L_{\l,k}\odot L_{\l, k}^\omega$ defined above
satisfy the assumptions of Lemma~\rom{10.1} and thus endow $\Ok$ with
the structure of a hermitian category.
\endproclaim

\proclaim{Lemma 10.6} Equivalence of categories $\Cal C(\g,
\varkappa)\simeq \Ok$ constructed in \cite{F} \rom{(}see
Theorem~\rom{9.6)} preserves the hermitian structure.
\endproclaim
\demo{Proof} The proof is based on the analysis of formulas defining the
equivalence of categories in \cite{KL3}.
\enddemo

\remark{Remark 10.7} So far, this definition depends on the
normalizations of the maps $\phi_\l$ which are defined up to a real
constant; however, as  in the quantum group case, the involution defined
above can be related with the usual compact involution on $\ghat$,
which is much more usual in the physical literature, and this allows
us to define the inner product up to a positive constant. We discuss this in
the Appendix. The natural conjecture, parallel to the quantum group
case, is that so defined inner product on the spaces of morphisms is
positive definite. So far, we have no proof of it except for $\sltwo$
case where it can be checked by direct calculation.
\endremark

\subhead Existence theorem \endsubhead

Now that we have defined the structure of a hermitian modular tensor
category on $\Ok$, the general theory developed in Section~2
immediately yields the existence of a modular invariant hermitian form
on the spaces of morphisms. 
 We formulate this result as a theorem. As
before, we denote

$$H=\bigoplus_{\l\in C}DL_{\l,k}\odot L_{\l, k}.
\tag 10.2$$ 

Let us introduce the following notation:
$$W_\l^{\mu\nu}=\Hom_{\Ok}(L_{\l,k}, L_{\mu,k}\odot L_{\nu,k}).
$$

\proclaim{Theorem 10.8} There exists a hermitian  form 
 on each of the spaces $W_\l^{\mu\nu}$ such that the resulting
form  on the space 
$$\Hom_{\Ok} (H, L_{\mu,k})=\bigoplus_{\l\in C} W_{\l^*\l}^\mu=
\bigoplus_{\l\in C} W_\l^{\l\mu}
\tag 10.3$$ 
is modular invariant.     
\endproclaim
\demo{Proof} Follows from the fact that $\Ok$ has a structure of
modular tensor category and constructions of Section~2.\qed\enddemo

 We will call this form ``the inner
product'', even though, as was noted before, we have no proof that
this form is positive definite.
Note that the  definition of this inner product depends on the choice of
the maps $\phi_\l:\1\to  L_{\l,k}\odot L_{\l, k}^\omega$ which were
used in the definition of hermitian structure. 
However, the inner product on the space $\Hom_{\Ok} (H, L_{\mu,k})$
in fact depends only on the choice of identification $\phi_\mu$
 and does not depend on the choice of $\phi_\l$
for each $\l$. Thus, the inner product on this space is defined
uniquely up to a constant factor. 

Recall (see Theorem~9.1) that the spaces
$W_\l^{\mu\nu}=\Hom_{\Ok}(L_{\l,k}, L_{\mu,k}\odot L_{\nu,k})$ are
identified with spaces of conformal blocks for a 3-punctured sphere;
thus, the construction above defines an inner product on the spaces of
conformal blocks  on a 3-punctured sphere.

To define an inner product on other Riemann surfaces with marked points,
we use the gluing, which allows us to represent the space of
conformal blocks assigned to a Riemann surface as a sum of tensor
products of the 3-point conformal blocks. Once we choose such a
representation, we define the norm by the rule $\|\Phi_1\o
\Phi_2\|= \|\Phi_1\|\cdot \|\Phi_2\|$. 

However, it is not clear why this inner product  is well-defined. Even for 
the conformal blocks on a sphere with 3 punctures, we have  so far defined
the inner product not for all 3-punctured spheres but only for some
contractible subspace $\Cal D$ in the moduli space. For other Riemann
surfaces, it is even more complicated, since there are different ways
to cut a surface in trinions. For example, there are different ways to
represent a torus with one puncture as a result of gluing together two
holes of a three-punctured sphere. Each of these ways gives an
identification of the space of conformal blocks with the space
$\bigoplus_{\l\in C} W_\l^{\l\mu}$ (see Lemma~9.4), and different
identifications give rise to the  action of the modular group on this
space -- cf. Theorem~9.5.

Now we come to the main result of our paper:

\proclaim{Theorem 10.9} There exists an inner product on
the spaces of conformal blocks on arbitrary Riemann surface 
such that with respect to this inner product, the natural flat
connection on the space of conformal blocks is
unitary. 
\endproclaim

\remark{Warning} As before, we do not claim that this inner product is
positive definite; thus, the untitarity should be understood
formally. 
\endremark

\demo{Proof} It follows from the general result of Moore and Seiberg
(see \cite{MS1, Section 5}) that it suffices to define such an inner
product on the sphere with $n$ punctures and on the torus with one
puncture. Moreover, to define the latter inner products it suffices
to check that the commutativity and
associativity morphisms, as well as the action of the modular group,  
are unitary with respect to the inner product defined on
$W_\l^{\mu\nu}$. But we have proved that in any modular tensor
category,  the inner product
defined in Section~2 satisfies these conditions (see Theorems~2.4, 2.5). 
\qed
\enddemo

This theorem is very important for Conformal Field Theory, since
existence of such an inner product is one of the axioms of CFT; thus,
Theorem~10.8 claims  that this axiom is satisfied in the WZW
model. To the best of my knowledge, this has not been proven in general
case so far. We refer the reader to the paper \cite{FG} for review of
known results and explicit constructions related to the  inner
products in higher genera.

Another important corollary of the equivalence of categories discussed
in Theorem~9.6 and results of Section~5 is the following theorem.

\proclaim{Theorem 10.10}
Let $\g=\sln$. Consider the category $\Ok$ with $\varkappa=K+kh\v$ for
some $K, k\in \Z_+$. 
Consider the space of conformal blocks on a torus with
one puncture to which the representation $L_{\mu,K+(k-1)h\v},
\mu=n(k-1)\omega_1$ is assigned \rom{(}here $\omega_1$ is the first
fundamental weight\rom{)}. Then the action of $SL_2(\Z)$ on this space
in some basis is given by  formulas \rom{(5.8)}; in particular, it
is written in terms of special values of Macdonald's polynomials at
roots of unity. 
\endproclaim

This theorem also gives the action of $SL_2(\Z)$ on the affine Jack
polynomials, introduced in \cite{EK}.  In this paper, for every $K,
k\in\Z_+, \l\in C_K=\{\l\in P^+| \<\l, \theta\v\>\le K\}$ we defined a
function $J_{\l,K}(h, u,\tau), h\in\h, u\in \C, \tau\in\C, Im\ \tau>0$
(see definition in \cite{EK, Section 8}) which we called the
(normalized) affine Jack polynomial; this name can be misleading,
since they are not polynomials but rather theta-functions, but it was
chosen since they are natural generalizations of Jack polynomials. We
proved that for $\g=\sln$, $J_{\l,K}$ can be calculated as suitable
renormalized traces of intertwining operators for the corresponding
affine Lie algebra, and that the space spanned by $J_{\l,K}$ (with
fixed $K,k$) is invariant under the action of $SL_2(\Z)$. Now, using
the previous results, we can calculate this action for $\g=\sln$.

\proclaim{Theorem 10.11} In the assumptions of the previous theorem,
let $J_{\l,K}$ be the normalized affine Jack polynomials as
in \cite{EK, Section 8}. Then 

$$\aligned
&J_{\l,K}(h,u,\tau+1)=q^{(\l+k\rho,\l+k\rho)-k\varkappa(\rho, \rho)/n},\\
&\tau^{-j}J_{\l,K}\bigl(\frac{h}{\tau},u-\frac{(h,h)}{2\tau},
			-\frac{1}{\tau}
		  \bigr)=
\sum_{\mu\in C_K}S_{\mu\l}J_{\mu,K}(h,u, \tau),
\endaligned
\tag 10.4 
$$
where $q=e^{\pi\i/\varkappa}, \varkappa=K+nk$, $S_{\l\mu}$ is given by formula \rom{(5.8)} in \cite{K}, and 
$j=-K(k-1)(n-1)/2\varkappa $.
\endproclaim

\demo{Proof} Since $J_{\l,K}$ are (up to a renormalization) the traces
of the intertwining operators for $\widehat{\sln}$, the action of the
modular group on these polynomials is up to simple renormalization the
same as the action on the space of intertwiners, or, equivalently, on
the space of conformal blocks. Combining it with the previous theorem,
we get formulas (10.4). The factor $\tau^{-j}$ appears in the formula
because the the tori $E_\tau=\C/\Z+\tau\Z$ is, of course, isomorphic
to $E_{-1/\tau}$, but the local parameters around zero, inherited from
$\C$, are different in these two realizations. The factor $\tau^{-j}$
accounts for the dependence of the space of conformal blocks on the
choice of local parameter.
\qed\enddemo

\head 11. Explicit formulas
\endhead

Our next goal is to give as explicit formulas as possible for this
inner product. First, let us introduce some notations.  

For every $\g$-homomorphism $g:V_1\to V_2\o V_3$ denote by $g^o:
V_2\to V_1\o V_3^*$ the image of $g^\omega$ (see (10.1)) under the
canonical isomorphism $\Hom_\g (V_1^*, V_3^*\o V_2^*)\simeq \Hom_\g
(V_2, V_1\o V_3^*)$ (as before, we assume that we have chosen
identifications $V^\omega\simeq V^*$). In a similar way, for every
$\Ok$-homomorphism $\Phi: V_1\to V_2\odot V_3$ we define $\Phi^o$ by
$(\Phi^g)^o =\Phi^{(g^o)}$ (see Example~8.7 for notations). 

Now we can  rewrite the definition of
$\|\Phi\|^2$ given in Section~2 in the following form:

\proclaim{Lemma 11.1} Let $\Phi\in \Hom_{\Ok} (V_1, V_2\odot
V_3)$. Then 

$$
\|\Phi\|^2=\frac{1}{(\dim_\varkappa V_1 \dim_\varkappa  V_2
\dim_\varkappa V_3)^{1/2}}\quad
\vcenter{\epsfbox{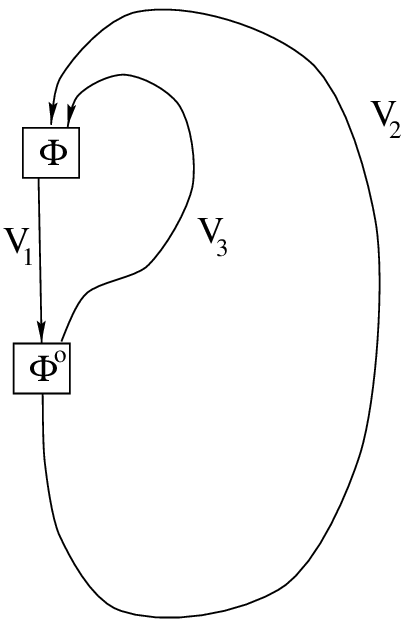}}.\tag 11.1$$

\endproclaim

\demo{Proof} Recalling the definition (2.2) of the inner product  and
Example~10.4, we see that it suffices to prove the following lemma: 

\proclaim{Lemma 11.2} 
$$\epsfbox{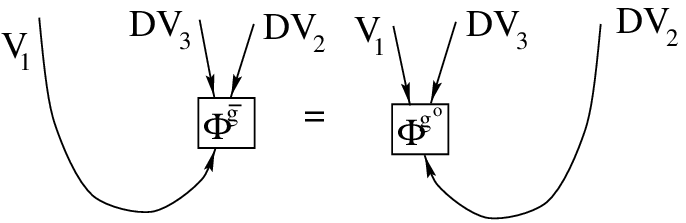} $$
\endproclaim

This Lemma follows from the fact that on the space $(V_1[0]\o V_2[0]\o
V_3[0])^\g$ the operators $\Omega_{12}$ and $\Omega_{23}$ commute and
thus, according to Corollary~9.2, the associativity isomorphism is
trivial. 
\qed\enddemo

Let us additionally assume that $V_1, V_2, V_3$ are irreducible:
$V_1=L_{\l, k}, V_2=L_{\mu,k}, V_3=L_{\nu,k}$. Denote by $\Psi$
composition of homomorphisms $\Phi, \Phi^o: \Psi=(\Phi\odot
1)\Phi^o\in \Hom_\Ok (V_2, (V_2\odot V_3)\odot DV_3)$. It follows from
the definition of the  associativity isomorphism in
$\Ok$ that $\Psi$ can be considered as a flat  section of the
bundle of conformal blocks $W(\Cal X,V_i)$, where $\Cal X$ is $\C P^1$
with the marked points $0, z_1, z_2, \infty$ and  representations
$V_2, V_3, DV_3, V_2$ assigned to these points respectively. This flat
section is uniquely defined by the following condition: if we denote
by $\<\Psi(z_1, z_2)\>$ the corresponding correlation function then it has the
following asymptotic as $z_1/z_2\to 0$: if $\Phi=\Phi^g$ then  

$$\<\Psi(z_1, z_2)\>\sim (g\o 1) g^o
z_1^{\Delta_\l-\Delta_\mu-\Delta_\nu} z_2^{\Delta_\mu-\Delta_\l-\Delta_\nu} $$

\proclaim{Theorem 11.3} In the notations above, we have the following
identity  of  $\g$-homo\-morphisms $L_\mu\to L_\mu\o L_\nu\o
L_\nu^*$: 

$$
\|\Phi\|^2 \Id_{L_\mu}\o i_\nu =
\biggl(\frac{\dim_q L_\mu\dim_q L_\nu}
      {\dim_q L_\l }
\biggr)^{1/2}
\lim_{z_1, z_2\to 1} (z_1-z_2)^{2\Delta_\nu}
\<\Psi(z_1, z_2)\>,
\tag 11.2, $$
where $i_\nu:\C\to  L_\nu\o L_\nu^* $ is the canonical map of
$\g$-modules.  
\endproclaim

\demo{Proof} 

It follows from the existence of the associativity isomorphism that $\Psi$ can
be rewritten as the sum 

$$\epsfbox{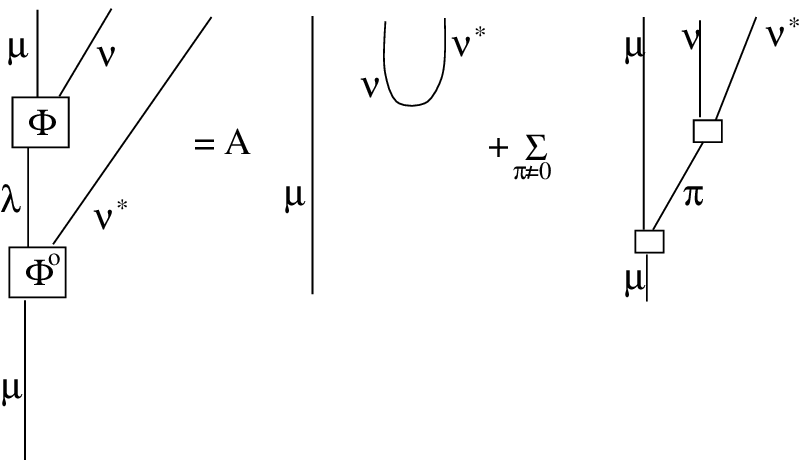}\tag 11.3$$ 
for some $A\in \C$ and some intertwiners in the boxes that are left
blank. Substituting this in the expression for $\|\Phi\|^2$ given in
Lemma~11.1, we see that 

$$\|\Phi\|^2=\biggl(\frac{\dim_q L_\mu\ \dim_q L_\nu} 
{\dim_q L_\l}\biggr)^{1/2} A.
\tag 11.4$$

As was mentioned before, coefficients in the right hand side of (11.3)
are related with asymptotic expansion of $\<\Psi\>$ as $z_1\to z_2$
(this is how the associativity morphism is defined). More precisely,
the first term in (11.3) gives the asymptotics $A
(z_1-z_2)^{-2\Delta_\nu}\Id_{L_\mu}\o i_\nu$, where $i_\nu: \C\to
L_\nu\o L_\nu^*$ is the canonical embedding and all other terms have
highest term of the asymptotics $(z_1-z_2)^{-2\Delta_\nu+\Delta_\pi}$.
Since $\Delta_\pi>0$ for $\pi\ne
0$, these terms give zero contribution to the limit, and therefore

$$\lim_{z_1, z_2\to 1}(z_1-z_2)^{2\Delta_\nu}\<\Psi(z_1, z_2)\>= A
\Id_{L_\mu}\o i_\nu.$$ 

\qed\enddemo

Thus, the calculation of the norm reduces to computation of the limit
(11.2). Since the asymptotics of $\<\Psi\>$ in the limit $|z_2|\gg
|z_1|$ is known, this is equivalent to calculation of Drinfeld's
associator. Note that this shows that the structure of $L_{\l,k}$ as a
module over $\ghat$ is not important here; in particular, this inner
product can be as well calculated for $\varkappa\notin \Q$, and the
answer is a meromorphic function of $\varkappa$, which may have poles
or zeroes only at rational points.  Since there are integral formulas
for the solutions of the KZ equations, the answer can be always
written in terms of asymptotics of some integrals. In the case of
$\g=\sltwo$ these asymptotics reduce to the Selberg integrals and are
given by certain products of $\Gamma$-functions (see \cite{V1}). For
arbitrary Lie algebras, the inner product can also be written in terms
of integrals of Selberg type, but in general can not be reduced to
gamma-functions.

\remark{Remark} In the case where all the spaces of conformal blocks
are zero or one dimensional (as happens for $\sltwo$), formula (11.2) (in
different form) for the inner product was suggested by Moore and Seiberg (see
\cite{MS2, Exercise 6.6}). \endremark

\head 12. Example: $\g=\sltwo$.\endhead

In this section we give explicit formulas for the norms of vertex
operators for $\sltwo$. As was noted in the end of the last section,
it suffices to calculate these norms for $\varkappa\notin \Q$. This
can be done using integral formulas for the solutions, given by
Schechtman and Varchenko \cite{SV}. Moreover, in this case the
asymptotics of the integrals can be calculated explicitly, using
the Selberg integral, and the answer is written in terms of
$\Gamma$-functions (see \cite{V1}). This can be used to find the
matrices which relate asymptotics of solutions in different asymptotic
zones.  Fortunately, this work has already been done by Varchenko in
\cite{V1}, in which there is a construction of equivalence of
categories $\Ok$ and $Rep\ \g$ for $\sltwo$.

In this section we only consider
$\g=\sltwo$. We identify the Cartan subalgebra with $\C$ so that
$\a\mapsto 1$, where $\alpha$  is the positive root; thus, $P\simeq
\frac{1}{2} \Z$. 

 Let us assume that we are given $K\in \C\setminus\Q$, and we
let $\varkappa=K+2$. Let $\Cal O_\varkappa$ be the category of
highest-weight modules of level $K$  over $\widehat{\sltwo}$ as in \cite{KL1}. 
It is known that this category is semisimple, and simple objects in
this category are precisely irreducible modules $L_{\l,K}$, which for
$\varkappa\notin \Q$ coincide with Weyl (induced) modules:
$L_{\l,k}=V_{\l, K}, \l\in P^+\simeq \frac{1}{2} \Z_+$. The same
constructions as before define on this category a structure of a
ribbon category, which is hermitian for $K\in \R$. Note that since we
assumed $\varkappa\notin \Q$, we have an isomorphism
$\Hom_{\Cal O_\varkappa}(V_1,V_2\odot V_3)\simeq \Hom_\sltwo(V_1[0], V_2[0]\o
V_3[0])$ (compare with (9.5)). Thus, this category could be described
without any reference to affine Lie algebra at all, just in terms of
finite-dimensional representations of $\sltwo$ and
Knizhnik-Zamolodchikov equation, as was done in \cite{Dr3, Dr4}.

For every triple $j, j_1, j_2\in \frac{1}{2} \Z_+$ such that $
|j_1-j_2|\le j\le j_1+j_2, j+j_1+j_2\in \Z$ there exists a unique up
to a constant $\sltwo$-morphism $g_j^{j_1 j_2}:V_j\to V_{j_1}\o
V_{j_2}$, and corresponding $\Cal O_{\varkappa}$-morphism $\Phi_j^{j_1
j_2} = \Phi^{g_j^{j_1 j_2}}$. 

Let us for simplicity consider the case $j=j_1, j_2=k$ for some
$j,k$. Then we can fix normalization of $g_j^{jk}$ by fixing some
vector of weight zero  $u\in V_k$ and requiring that $g_j^{jk} v_j
=v_j\o u +\dots$, where $v_j$ is a highest weight vector in $V_j$. 

The main result of this section, which can be obtained by the use of
explicit formulas in \cite{V1}, is the following theorem:

\proclaim{Theorem 12.1} Let $\varkappa\in \R\setminus \Q$. Then: 

$$\|\Phi_j^{j k}\|^2=c(k,\varkappa)
 \prod_{i=1}^{k} 
	\frac{[2j+1+i]}
	     {[2j+1-i]}
	\frac{\[2j+1+i\]^2}
	     {\[2j+1-i\]^2}
\tag 12.7$$
where the constant $c(k,\varkappa)$ does not depend on $j$ and 

$$\[x\]=\Gamma\biggl(-\frac{x}{\varkappa}+1\biggr).\tag 12.8$$
\endproclaim

Since $\[x\]$ is well-defined and non-zero for $x\ne n\varkappa,
n=1,2,..$ we see that for $\varkappa\notin \Q$ the norm $\|\Phi_j^{j
k}\|$ is well-defined and non-zero. One can also easily check that if
$K\in \Z_+, j,k\le K/2$ (this last condition ensures that $L_{j, K}, L_{k,K}$
are integrable) then  $\|\Phi_j^{j k}\|$ is well-defined and non-zero
iff $k\le 2j\le K-k$, which is exactly the condition for the space of
intertwiners $\Hom_{\Ok}(L_{j,K}, L_{j,K}\odot L_{k,K})$ to be
non-zero.

Note also that as $\varkappa\to \infty$, formula (12.7) coincides (up to a
constant) with the Macdonald's inner product identity (5.7) for
$\sltwo$ (with $q=1$).

This theorem justifies the following conjecture, which is a natural
generalization of the Macdonald's inner product identities to affine
root systems.  Let $\g=\sln, k\in \Z_+, \varkappa\in
\R\setminus \Q$. For $\l\in P_+$, choose $\hat\Phi_\l\in \Hom_{\Cal
O_\varkappa} (V, V\odot U)$, where $V=V_{\l+(k-1)\rho, \varkappa-h\v},
U=V_{n(k-1)\omega_1, \varkappa-h\v}$ so that $\hat\Phi(v)=v\o u_0+\dots$,
where $v$ is a highest-weight vector in $V$ and $u_0\in U[0]$ is some
fixed zero-weight vector. This is the natural affine analogue of the
intertwiners $\Phi_\l$ used in the quantum group case to obtain
Macdonald's polynomials (see \cite{K, Section 5} and references
therein); as before, one can show that $\Hom_{\Cal O_\varkappa} (V,
V\odot U)$ is one-dimensional and thus the condition $\hat\Phi_\l(v)=v\o
u_0+\dots$ uniquely defines $\hat\Phi_\l$. The intertwiners $\hat \Phi_\l$
were used in \cite{EK} to construct an affine analogue of Jack's
polynomials. 

As before, $\Cal O_\varkappa$ has a natural structure of a hermitian
category and thus we have an inner product on the spaces of
morphisms. 

\proclaim{Conjecture 12.2 \rm (Affine inner product identities)} 
Let $\g=\sln, \varkappa\in \R\setminus \Q$, and let $\Phi_\l:V\to
V\odot U$ be as above. Then: 

$$\|\hat\Phi_\l\|^2=c(k,\varkappa)\prod_{\a\in R^+}
		\frac{[(\a, \l+k\rho)+i]}
		     {[(\a, \l+k\rho)-i]}
		\left(\frac{\[(\a, \l+k\rho)+i\]}
		     		{\[(\a, \l+k\rho)-i\]}
		\right)^2,
\tag 12.9
$$
where the constant $c(k,\varkappa)$ does not depend on $\l$. 
\endproclaim

In the case $\g=\sltwo$ this conjecture coincides with the statement
of Theorem~12.1. 

This formula is indeed a natural analogue of Macdonald's inner product
identity (5.7) for the following reason: if we write weights for
$\ghat$ in the form $\hat\l=\l+k\Lambda_0$ where $\Lambda_0$ is dual
to the central element $c$, and denote $\hat\rho=\rho +h\v\Lambda_0$
then the right hand side of (12.9) is exactly the regularization of
the meaningless expression

$$\prod_{\a\in \widehat R^+_{re}}\prod_{i=1}^{k-1} 
\frac{(\a, \hat\l+k\hat\rho)+i}
	{(\a, \hat\l+k\hat\rho)-i}
$$ 
obtained by replacing products of the form 
$\prod_{i=0}^\infty \frac{1}{x+i}$ with $\Gamma(x)$. 

This conjecture is closely related to the expressions given by
Gaw{\c e}dzki et al in \cite{G,FG, FGK}. In these papers they suggest a
construction of the inner product on the spaces of conformal blocks
based on the fact that these spaces are the state spaces of the quantum
Chern-Simons theory. The inner product is obtained by regularization
of certain infinite-dimensional integrals, and the final answer is
given by a finite-dimensional integral of the following type

$$\int_{\C^k}|\omega|^2,$$
where $\omega$ is the  differential form 
which appears in the integral formulas of Schechtman and Varchenko
(up to $\varkappa\mapsto -\varkappa$). 

There are reasons to believe that in the case $\g=\sltwo$ integrals of
these type can be calculated explicitly (see \cite{A} for a simplest
example of computation of this type), and the answer can be expressed
in terms of $\Gamma$-functions similarly to the Selberg integral. We
expect that this procedure would yield the same answer as the one
given by formula (12.7) above.

\head Appendix\endhead

Here we give another description of the hermitian structure on $\Ok$,
based on compact involution of $\g$. This construction is parallel to
the one for $\Ug$ (see Section~7).

Let $\omega_c$ be the compact involution on $\g$, i.e. an antilinear
Lie algebra automorphism given on the generators  by 
$\omega_c(e_i)=-f_i,
\omega_c(f_i)=-e_i, \omega_c(h)=-h$. Let us extend it to an involution
on $\ghat$ by letting $\omega_c(x[n])=\omega_c(x)[n],\omega_c(c)=c$
({\bf Warning: this is different from the standard compact involution
on $\ghat$!}). We define $V^{\omega_c}, \Phi^{\omega_c}$ similar to
above. Note that for any $\l$, $L_{\l,k}^{\omega_c}\simeq L_{\l^*,
k}$, and thus, from complete reducibility, we have $V^{\omega_c}\simeq
DV$ for any $V\in \Ok$, though these isomorphisms are not canonical.

For any $\Cal X$ as in the beginning of Section~8, let $\overline{\Cal
X}$ be the following collection of data:

\roster 
\item The curve $\overline{X}$ which is the same curve as $X$ but with
the opposite complex structure (thus, if $f$ is a meromorphic function on
$X$ then $\bar f$ is a meromorphic function on $\overline{X}$).

\item The marked points $z_i$ are the same as for $\Cal X$.

\item The new local parameters  are $\overline{w_i}$ (here we consider
$w_i$ as holomorphic function on $X$ vanishing at $z_i$). 
\endroster

\example{Example} Let $\Cal X$ be the Riemann sphere (8.8). Then
$\overline{\Cal X}$ is isomorphic to $\C P^1$ with global parameter
$w$, marked points $\overline{z_i}, \infty $ and local parameters
$\bar z_i-w, 1/w$ 
respectively. 
\endexample

Assume that we are given  $\Cal X$  as before and we have
representations $V_i$ assigned to the marked points. Let $\Phi\in
W(\Cal X; V_1,V_2,\dots, V_n)$ be a conformal block, i.e. a map 

$$\o_{In}V_i\to \hat\o_{Out}V_i$$
which commutes with the action of the algebra $\Gamma$ of meromorphic
$\g$-valued functions on $\Cal X$ (see Definition~8.1).

 Let us denote
by $\Phi^{\omega_c}$ the same map considered as a map 

$$\o_{In}V_i^{\omega_c}\to \o_{Out}V_i^{\omega_c}.$$

\proclaim{Lemma A.1} The mapping $\Phi\mapsto \Phi^{\omega_c}$ is an
involutive isomorphism\newline
 $W(\Cal X; V_1, \dots, V_n)\to W(\overline{\Cal
X}; V_1^{\omega_c},\dots,V_n^{\omega_c})$.
\endproclaim
\demo{Proof} This is obvious from the fact that $f(z)\mapsto
\omega_c(f(z))$ identifies the algebras $\Gamma$ of meromorphic
functions on $\Cal X$ and $\overline{\Cal X}$. 
\qed\enddemo
 
In particular, if $\Cal X$ is the Riemann sphere (8.8) then
$\Phi^{\omega_c}$ will be a conformal block on the sphere with marked
points $\bar z_1,\dots, \bar z_n$. This shows that we have a natural
isomorphism $(V\odot W)^{\omega_c}\simeq V^{\omega_c}\odot W^{\omega_c}$.

Now we can define the hermitian structure on $\Ok$ similarly to the
definitions before (see Lemma 10.1); the only non-trivial part is that
we define isomorphisms $(V\odot W)^{\omega_c}\simeq W^{\omega_c}\odot
V^{\omega_c}$ to be given by the composition of the natural
isomorphism $(V\odot W)^{\omega_c}\simeq V^{\omega_c}\odot
W^{\omega_c}$ and the commutativity morphism $Pe^{\pi\i
\Omega/\varkappa}$, where $P$ is the transposition. Also, we define
the maps $\phi_\l:\1\to L_{\l,k}\odot L_{\l,k}^{\omega_c}$ by
$\phi_\l(1)=e^{\pi\i\Delta_\l}v_{\l,k}\o v_{\l,k}^{\omega_c}+\dots$,
where $\Delta_\l$ is defined by (9.3).

\proclaim{Lemma A.2} So defined functor $\omega_c$ satisfies all the
conditions of Lemma~\rom{10.1} and thus gives rise to a hermitian
structure on $\Ok$.
\endproclaim

In fact, it turns out that this hermitian strucutre is equivalent to
the one described in Section 10.  As in the classical case, this
relies on the existence of an analogue of the longest element of Weyl
group.

For every module $V\in \Ok$ define the map 

$$\Omega=w_0e^{-\pi\i L_0} :V\to V,\tag 10.2$$
where $w_0$ is the longest element of Weyl group for $\g$, which acts
on every finite-dimensional representation of $\g$ (see, e.g.,
\cite{KL3, 19.4}); thus, it also acts on every module in $\Ok$.  Its
$q$-analogue (which we denoted also by $\Omega$) was discussed in
 Section~7.

One can easily check that 

$$\gathered
\omega(x)=\Omega^{-1}\omega_c( x)\Omega
	=\omega_c(\Omega^{-1} x\Omega),\quad x\in \ghat\\
\Omega^2 = Z e^{-2\pi\i L_0}
\endgathered
\tag 10.3
$$
where $Z=w_0^2$ was discussed in Section~7. Recall that $Z$ is a
central element which  acts by $\pm 1$  in any irreducible $\g$-module
and satisfies $Z|_{V\o W}=Z_V\o Z_W$; thus, it also acts by constant
in any simple $\ghat$-module from category $\Ok$ and satisfies
$Z_{V\odot W}=Z_V\odot Z_W, Z^2=1$. 

Thus, for every representation $V\in \Ok$ the map $\Omega: V\to V$,
considered as a map $V^{\omega}\to V^{\omega_c}$, is an isomorphism of
$\ghat$-modules; similarly, we can identify
$\overline{\Phi}=\Omega^{-1} \Phi^{\omega_c}\Omega$.

\enddocument